\documentclass[11pt]{article}
\usepackage{amsbsy}
\usepackage{amssymb}
\newcommand{\be}{\begin{equation}}
\newcommand{\ee}{\end{equation}}
\newcommand{\bdm}{\begin{displaymath}}
\newcommand{\edm}{\end{displaymath}}
\newcommand{\beqr}{\begin{eqnarray}}
\newcommand{\eeqr}{\end{eqnarray}}
\newcommand{\beqrn}{\begin{eqnarray}}
\newcommand{\eeqrn}{\end{eqnarray}}

\def\t{\theta}
\def\x{\chi}

\begin{document}
\title{Euler's triangle and the decomposition of tensor powers of adjoint representation 
of $A_1$ Lie algebra}
\author{A.M. Perelomov}
\date{}
\maketitle
\begin{center}
${\it  Institute\, for\, Theoretical\, and\, Experimental\, Physics, 117259, Moscow, Russia.}$
\end{center}
\begin{abstract}
\noindent
We consider the relation between Euler's trinomial problem and the problem of decomposition
of tensor powers of adjoint representation of $A_1$ Lie algebra. 
By using this approach, some new results for both problems are obtained.
\end{abstract}

\noindent
{\bf 1. Introduction.} 
\medskip 

In 1765 Euler [Eu 1765/1767] investigated the coefficients of trinomial
\be
(1+x+x^2)^n=\sum _{k=-n}^n {a_n^{(k)}x^{n+k}}.
\ee 
For central trinomial coefficients $a_n^{(0)}$ he found the generating function and two term recurrence relation.  
For discussion of properties $a_{n}^{(k)}$ see [Ri 1974].

Let us change $x$ to $\exp(i\t)$ and rewrite the left hand side of (1) as
\be
(1+x+x^2)^n=x^n\,X^n,\quad X=1+2\,\cos\t\,.
\ee

Note that $X$ is character of adjoint representation of Lie algebra $A_1$. 
So, Euler's problem is equivalent to the problem of multiplicities of weights of representation $X^n\,.$
We consider also related to this problem
of decomposition of $X^n$ into irreducible representations of Lie algebra $A_1\,.$

\medskip\noindent
{\bf 2. Euler's triangle.}

It is evident that $a_n^{(-k)}=a_n^{(k)}$. So, it is sufficient to consider only quantities 
$a_{n}^{(k)}$, $k\geq 0$. It is convenient to arrange these coefficients in triangle.
We give here the table of these numbers till $n=10$:

\medskip
\begin{tabular}{lllllllllllll}
   {\bf n/k}&    &            0&     1&     2&     3&     4&     5&    6&    7&   8&   9&  10   \\
           0&    &      {\bf 1}&      &      &      &      &      &     &     &    &    &      \\
           1&    &      {\bf 1}&     1&      &      &      &      &     &     &    &    &      \\
           2&    &      {\bf 3}&     2&     1&      &      &      &     &     &    &    &       \\
           3&    &      {\bf 7}&     6&     3&     1&      &      &     &     &    &    &      \\
           4&    &     {\bf 19}&    16&    10&     4&     1&      &     &     &    &    &      \\
           5&    &     {\bf 51}&    45&    30&    15&     5&     1&     &     &    &    &      \\
           6&    &    {\bf 141}&   126&    90&    50&    21&     6&    1&     &    &    &       \\ 
           7&    &    {\bf 393}&   357&   266&   161&    77&    28&    7&    1&    &    &       \\
           8&    &   {\bf 1107}&  1016&   784&   504&   266&   112&   36&    8&   1&    &       \\
           9&    &   {\bf 3139}&  2907&  2304&  1554&   882&   414&  156&   45&   9&   1&       \\
          10&    &   {\bf 8953}&  8350&  6765&  4740&  2850&  1452&  615&  210&  55&  10&   1       
\end{tabular}

From (1) it follows immediately the three term recurrence relation
\be
a_{n+1}^{(k)} =a_n^{(k-1)}+a_n^{(k)}+a_n^{(k+1)}.
\ee

Let us introduce the generating 
function $F(t)$ for central trinomial coefficient $a_n=a_n^{(0)}$
\be
F(t)=\sum_{n=0}^{\infty } a_n\,t^n .
\ee
{\bf Theorem (Euler 1765)}. {\em The generating function $F(t)$ has the form}
\be 
F(t)=\frac {1}{\sqrt{(1-2t-3t^2)}}
\ee
{\em and we have two term recurrence relation for coefficients} $a_{n}$
\be 
n\,a_{n}=(2n-1)\,a_{n-1}+3(n-1)\,a_{n-2}. 
\ee

We give here a very short proof of the first statement different from Euler's one.

\noindent
{\bf Proof.} Note that
\be 
a_n=\frac{1}{\pi} \int_0^{\pi}X^n\,d\t,\qquad X=(1+2\cos\t).
\ee
So,
\be
F(t)=\frac{1}{\pi} \int_0^{\pi} \frac{d\t}{1-t-2t\,\cos\t}. 
\ee
Evaluating this integral we obtain formula (5).

The second statement is a special subcase of more general statement.

\noindent
{\bf Theorem 1.} {\em We have  two term recurrence relation for coefficients} $a_{n}^{(k)}$
\be 
(n^2-k^2)\,a_n^{(k)} =n(2n-1)\,a_{n-1}^{(k)}+3n(n-1)\,a_{n-2}^{(k)}. 
\ee


\noindent
{\bf Proof.} We have 
\be
 a_n^{(k)}=\frac{1}{\pi } \int _0^{\pi } X^n\,\cos k\t\,d\t ,
\ee
and 
\be \int_0^{\pi } X^n\,\left[ \left( \frac{d^2}{d\t ^2} +k^2\right)\,\cos k\t \right]d\t =0
=\int_0^{\pi }\,\cos k\t \,\left[ \left( \frac{d^2}{d\t ^2}+k^2\right)X^n\right] d\t .
\ee
But,
\be 
\frac{d^2X^n}{d\t ^2}=-\,n^2X^n+n(2n-1)\,X^{n-1}+3n(n-1)\,X^{n-2}.
\ee
From this, it follows equation (9).

\medskip\noindent
{\bf Theorem 2.} {\em We have another two term recurrence relations for coefficients} $a_{n}^{(k)}$
\be 
(n+1)\,\left( a_n^{(k-1)}-a_n^{(k+1)}\right) =k\,a_{n+1}^{(k)}.
\ee
\be 
(n-k+1)\,a_{n}^{(k-1)}=ka_{n}^{(k)}+(n+k+1)\,a_{n}^{(k+1)},
\ee
\be 
(n-k+1)\,a_{n+1}^{(k)}=(n+1)(a_{n}^{(k)} +2a_{n}^{(k+1)}),
\ee
\be 
(n+k+1)\,a_{n+1}^{(k)}=(n+1)\,(a_{n}^{(k)}+2a_{n}^{(k-1)}).
\ee

\medskip\noindent
{\bf Proof.} From the identity
\be \int _0^{\pi } \left[ \frac{d}{d\t }\,(X^n\,\sin k\t )\right] d\t =0, 
\ee
we obtain equation (13). Combining this relation  with (3), we obtain equations (14)--(16).

\medskip\noindent
Note that from (3) it follows
$$
a_n^{(1)}=\frac12\,(a_{n+1}-a_n), \quad a_n^{(2)}=\frac12\,(a_{n+2}-2 a_{n+1}-a_n),
$$
\be
a_n^{(3)}=\frac12\,(a_{n+3}-3 a_{n+2} +2 a_n), \quad a_n^{(4)}=\frac12\,(a_{n+4}-4 a_{n+3}+ 2 a_{n+2}+ 4 a_{n+1}-a_n)\, .
\ee

\medskip\noindent
{\bf Corollary 1.}\,
Explicit expressions for quantities $a_n^{(n-k)}$ for small $k$\,
may be obtained from (9) and (14) and we have
\be
a_n^{(n-k)}=\frac{1}{k!} Q_k(n),
\ee
where $Q_k(n)$ is polynomial of degree $k$ in $n$. 

The recurrence relation for these polynomials follows from (14)
and we give the explicit expression for first ten polynomials.
\be
Q_{k+1}(n)=(n-k)Q_{k}(n)+k(2n-k+1)Q_{k-1}(n).
\ee

\be
Q_0=1;\qquad Q_1=n;\qquad Q_2=n(n+1); 
Q_3=(n-1)n(n+4);
\ee
$$ 
Q_4=(n-1)n(n^2+7n-6); 
$$
$$
Q_5=(n-2)(n-1)n(n+1)(n+12);
$$
$$
Q_6=(n-2)(n-1)n\,(n^3+18n^2+17n-120);
$$
$$
Q_7=(n-3)(n-2)(n-1)n\,(n^3+27n^2+116n-120);
$$
$$
Q_8=\,(n-3)(n-2)(n-1)n\,(n+1)(n+10)(n^2+23n-84);
$$
$$ 
Q_9=\,n(n-1)(n-2)(n-3)(n-4)\,(n^4+46n^3+467n^2+86n-3360);
$$
$$
Q_{10}=\,n(n-1)(n-2)(n-3)(n-4)\,(n^5+55n^4+665n^3-895n^2-16626n+15120).
$$

\medskip\noindent
{\bf 3. Decomposition of $X^n$ into irreducible representations.}

This problem is equivalent to expansion of $X^n$ in terms of characters of Lie algebra $A_1$.
\be
X^n=\sum _{k=0}^n b_n^{(k)} \x_k (\t) \,.
\ee
These characters are well known (see for example [We 1939].)
\be
\x _k =1 +2\cos (\t) + 2\cos (2\t) + ... + 2\cos (k\t)\,.
\ee
They are orthogonal
\be
\frac{1}{\pi} \int_0^{\pi}\x _k (\t)\,\x _l(\t)\,(1- \cos(\t))\, d\t  = \delta _{k,l}\,,
\ee 
and we have
\be
b_n^{(k)}=\frac{1}{\pi}\int_0^{\pi} X^n f_k(\t) d\t\,,\quad f_k(\t)=\cos(k\t)-\cos((k+1)\t)\,.
\ee
From this it follows the basic relation
\be
b_n^{(k)} = a_n^{(k)} - a_n^{(k+1)}\,,
\ee
the three term recurrence relation similar to relation (3)
\be
b_{n+1}^{(k)} =b_n^{(k-1)}+b_n^{(k)}+b_n^{(k+1)}\,,\quad {\rm for}\, n\geq 2\,, k\geq 1\,,
\ee
and other relations
$$
b_{n}=b_{n}^{(0)}=\frac{1}{2}(3 a_n - a_{n+1})\,, b_n^{(1)}=b_{n+1},\, b_{n}^{(2)}=b_{n+2}-b_{n+1} -b_n\,,
$$
\be
b_{n}^{(3)}=b_{n+3}- 2 b_{n+2} -b_{n+1}+b_n,\, b_{n}^{(4)}=b_{n+4}-3 b_{n+3} + 3 b_{n-1}\,.
\ee

\medskip\noindent
Here we have the triangle

\begin{tabular}{lllllllllllll}
   {\bf n/k}&    &            0&     1&     2&     3&     4&     5&    6&    7&   8&   9& 10  \\
           0&    &      {\bf 1}&      &      &      &      &      &     &     &    &    &      \\
           1&    &      {\bf 0}&     1&      &      &      &      &     &     &    &    &      \\
           2&    &      {\bf 1}&     1&     1&      &      &      &     &     &    &    &       \\
           3&    &      {\bf 1}&     3&     2&     1&      &      &     &     &    &    &      \\
           4&    &      {\bf 3}&     6&     6&     3&     1&      &     &     &    &    &      \\
           5&    &      {\bf 6}&    15&    15&    10&     4&     1&     &     &    &    &      \\
           6&    &     {\bf 15}&    36&    40&    29&    15&     5&    1&     &    &    &       \\ 
           7&    &     {\bf 36}&    91&   105&    84&    49&    21&    6&    1&    &    &       \\
           8&    &     {\bf 91}&   232&   280&   238&   154&    76&   28&    7&   1&    &       \\
           9&    &    {\bf 232}&   603&   750&   672&   468&   258&  111&   36&   8&   1&       \\
          10&    &    {\bf 603}&  1585&  2025&  1890&  1398&   837&  405&  155&  45&   9&   1       
\end{tabular}

\medskip\noindent
{\bf Theorem 3.} {\em The generating function $G(t)=\sum _{n=0}^{\infty} b_n t^n$ here has the form.}
\be 
G(t)=\frac{1}{2 t}\left(1-\frac {\sqrt{1-3t}}{\sqrt{(1+t)}}\right)\,.
\ee

\medskip\noindent
{\bf Proof.} Take into account identity
\be
\frac{1- \cos (\theta)}{1-t-2 t \cos(\theta)} = \frac{1}{2 t} (1-\frac{1}{1-t-2 t \cos (\theta)})\,.
\ee
Then the proof is reduced to the proof for  $F(t)$\,.
We have also the recurrence relation
that follows from (6) and $b_n=a_n -a_{n}^{(1)}$
\be
(n+1)b_n=(n-1)(2 b_{n-1}+3 b_{n-2})\,.
\ee

\medskip\noindent
{\bf Theorem 4.}

\noindent
{\em We have here four term recurrence relation}
\be
A_{nk} b_{n}^{(k)}+ B_{nk} b_{n-1}^{(k)}+ C_{nk} b_{n-2}^{(k)}+ D_{nk} b_{n-3}^{(k)}+ E_{nk} b_{n-4}^{(k)} =0\,, 
\ee
where
\be A_{n,k}=(n^2-(k+1)^2)(n^2-k^2); B_{n,k}=-2n(2n-1)(n+k)(n-k-1);
\ee
\be C_{n,k}=-2n(n-1)(n^2-2n+3-3k(k+1));
\ee
\be D_{n,k}=6n(n-1)(n-2)(2n-3); E_{n,k}=9n(n-1)(n-2)(n-3).
\ee

\noindent
{\bf Proof.}
We have 
\be 
b_n^{(k)}=\frac{1}{\pi} \int_0^{\pi}X^n f_k (\t) d\t
\ee
where
\be
X=(1+2 \cos(\t)), f_k(\t)=\cos(k\t)-\cos((k+1)\t)
\ee
and
\be A_k f_k(\t)=0, A_k=(\frac{d^2}{{d\t}^2}+k^2)(\frac{d^2}{{d\t}^2}+(k+1)^2)
\ee 
Integrating by part in (36) we get
\be
\frac{1}{\pi} \int_0^{\pi} f_k (\t) (A_k X^n) d\t =0\,\, {\rm and}\, (32) -- (35).
\ee

\medskip\noindent
{\bf Theorem 5.}

\medskip\noindent
{\em We have here three term recurrence relation}
\be (k+1)(n+1-k)b_n^{(k-1)}=(k(k+1)-n-1)b_n^{(k)}+k(n+k+2)b_n^{(k+1)}.\ee

\noindent
{\bf Proof.} This folows from (14) and the relation $b_{n}^{(k)}=a_{n}^{(k)}-a_{n}^{(k+1)}$

\end{document}